\documentclass{amsart}

\usepackage[T1]{fontenc}
\usepackage{lmodern}
\usepackage{amsmath,amssymb,amsthm,mathtools}
\usepackage{mathrsfs}
\usepackage[hidelinks]{hyperref}
\usepackage[nameinlink,capitalize]{cleveref}

\newcommand{\Tr}{\operatorname{Tr}}

\newcommand{\eps}{\varepsilon}

\newtheorem{theorem}{Theorem}[section]
\newtheorem{proposition}[theorem]{Proposition}
\newtheorem{lemma}[theorem]{Lemma}

\theoremstyle{definition}

\theoremstyle{remark}
\newtheorem{remark}[theorem]{Remark}

\title[Embeddings into Lorentz--Schatten ideals]
{Isomorphic embeddings of \(L_q(0,1)\) into
Lorentz--Schatten ideals}
\author[C. Gao]{Chenhao Gao}
\address{Institute for Advanced Study in Mathematics of HIT, Harbin Institute of Technology, Harbin, 150001, China}
\author[J. Huang]{Jinghao Huang}
\address{Institute for Advanced Study in Mathematics of HIT, Harbin Institute of Technology, Harbin, 150001, China}
\email{jinghao.huang@hit.edu.cn}
\thanks{This work was supported by the National Natural Science Foundation of China (NNSFC) under grants 12301160, 12471134, and 12671159.}
\date{}

\begin{document}

\begin{abstract}
Let \(0<p,q<\infty\).  We prove
\[
 L_q(0,1)\hookrightarrow\mathcal C_{p,q}
 \quad\Longleftrightarrow\quad q=2.
\]
As a consequence,
we obtain, for all \(0<p_1,q_1,p_2,q_2<\infty\),
\[
 L_{p_1,q_1}(0,1)\hookrightarrow\mathcal C_{p_2,q_2}
 \quad\Longleftrightarrow\quad p_1=q_1=2.
\]
\end{abstract}

\maketitle

\section{Introduction}

For \(0<p,q<\infty\), let \(\ell_{p,q}\) be the Lorentz sequence
space and let \(\mathcal C_{p,q}\) be its associated ideal of compact
operators on \(\ell_2\).  In particular, \(\mathcal C_{p,p}\) is the
Schatten--von Neumann class \(\mathcal C_p\).  Throughout the paper,
\(X\hookrightarrow Y\) means that \(X\) is isomorphic to a closed
subspace of \(Y\).

Arazy and Lindenstrauss \cite[Theorem~6]{AL} proved that
\(L_p(0,1)\) does not embed isomorphically into \(\mathcal C_p\) for
\(2<p<\infty\).  Their argument was extended by Astashkin, Huang,
and Sukochev \cite{AHS}, who established non-embedding results for
symmetric function spaces and operator ideals.  In the Lorentz case,
\cite[Theorem~5.1]{AHS} implies that
\(L_q(0,1)\not\hookrightarrow\mathcal C_{p,q}\) whenever
\(2<p\leq q<\infty\).

For \(q>2\), Jinghao Huang, Fedor Sukochev, and Dmitriy Zanin
\cite{HSZ} left open whether \(L_q(0,1)\) embeds isomorphically
into \(\mathcal C_{p,q}\) when \(0<p\leq2\) or \(p>q\).
We determine all pairs \((p,q)\) for which the embedding exists.

\begin{theorem}\label{thm:classification}
Let \(0<p,q<\infty\).  Then
\[
 L_q(0,1)\hookrightarrow\mathcal C_{p,q}
 \quad\Longleftrightarrow\quad q=2.
\]
\end{theorem}

For \(q>2\), \cref{thm:classification} extends the non-embedding
results of Arazy--Lindenstrauss and Astashkin--Huang--Sukochev to
\(p>q>2\) and \(0<p\leq2<q\).  We use Arazy's subsequence
decomposition theorem in the quasi-Banach setting of Huang,
Sukochev, and Yu \cite{HSY}.
For \(0<q<2\), the non-embedding follows directly from
\cite[Lemma~3.7]{HSZ} with \(p_1=q_1=q_2=q\) and \(p_2=p\),
since \(L_{q,q}(0,1)=L_q(0,1)\).

Section~\ref{sec:decomposition} recalls the subsequence decomposition
used in the proof for \(q>2\).  Section~\ref{sec:estimates}
establishes lower estimates for one-sided disjoint families and a
compression estimate for Hilbertian sequences.
In Section~\ref{sec:main-proof} we apply these estimates to a
perturbation of the Haar system and complete the proof of
\cref{thm:classification}.

\section{Subsequence decomposition}\label{sec:decomposition}

Let \(H=\ell_2\), with its standard basis \((e_j)\).  For a closed
subspace \(K\subset H\), write \(p_K\) for the orthogonal projection
onto \(K\).  We use the coordinate projections and square compressions
\[
 \pi_0=0,\qquad \pi_d=p_{[e_1,\ldots,e_d]}\quad(d\geq1),\qquad
 \mathscr R_d(x)=\pi_dx\pi_d\quad(d\geq0).
\]
For \(x\in\mathcal C_{p,q}\), we use the quasi-norm
\[
 \|x\|_{p,q}
 =
 \left(\sum_{n=1}^{\infty}
 n^{q/p-1}\mu_n(x)^q\right)^{1/q},
\]
where \((\mu_n(x))\) denotes the decreasing sequence of singular
values of \(x\).
Left and right multiplication by contractions are contractions on
every symmetric operator ideal.  Since \(q<\infty\), finite-rank
operators are dense in \(\mathcal C_{p,q}\); consequently
\begin{equation}\label{eq:square-compression}
 \|\mathscr R_d(x)-x\|_{p,q}\longrightarrow0\quad(d\to\infty)
 \qquad(x\in\mathcal C_{p,q}).
\end{equation}

A family \((x_i)\) is called left-disjoint (respectively,
right-disjoint) if its left (respectively, right) support
projections are pairwise orthogonal.
Unless stated otherwise, constants may depend on \(p,q\) and on the
embedding constants, but are independent of the dyadic level and
the number of steps in the construction.

The following form of Arazy's decomposition follows from
\cite[Theorem A.6 and Corollaries A.7--A.9]{HSY}, together with
\cite[Theorem 3.8]{HSY}.

\begin{proposition}\label{prop:hsy}
Let \(E\) be a separable symmetric quasi-Banach sequence space with
\(c_0\not\hookrightarrow E\), and let \((v_i)\subset\mathcal C_E\)
be a seminormalized \(\sigma\)-weakly null basic sequence.  Given
\(\eta>0\), there exist a subsequence, again denoted by \((v_i)\),
two families of finite-dimensional coordinate blocks
\((K_i')_{i\geq1}\) and \((L_i')_{i\geq1}\), and operators
\(a_i,b_i,c_i\) for \(i\geq2\), with the following properties.
The blocks within each family are pairwise orthogonal.  The
approximants \(z_i=a_i+b_i+c_i\) satisfy
\begin{align}
 a_i&\in p_{K_i'}\mathcal C_Ep_{L_1'},&
 b_i&\in p_{K_1'}\mathcal C_Ep_{L_i'},&
 c_i&\in p_{K_i'}\mathcal C_Ep_{L_i'},\label{eq:hsy-support}\\
 \left\|\sum_i t_i(v_i-z_i)\right\|_{\mathcal C_E}
 &\leq
 \eta\left\|\sum_i t_iv_i\right\|_{\mathcal C_E}
 \label{eq:hsy-relative}
\end{align}
for every finitely supported scalar sequence \((t_i)\).
Here \(K_1'\) and \(L_1'\) are the fixed blocks, and the chosen
subsequence and the sums are indexed by \(i\geq2\).  Moreover,
there is a number \(\beta\geq0\) such that
\begin{equation}\label{eq:hsy-hilbert}
 \left\|\sum_i t_i(a_i+b_i)\right\|_{\mathcal C_E}
 =\beta\|(t_i)\|_2 .
\end{equation}
\end{proposition}

\begin{proof}
Apply \cite[Corollary A.9]{HSY} with \(K_j=L_j=[e_j]\).
Passing to a subsequence, we obtain disjoint coordinate blocks
\(K_i',L_i'\) and operators
\[
 z_i=a_i+b_i+c_i
\]
which satisfy \eqref{eq:hsy-support} and \eqref{eq:hsy-relative}.
The representations in parts (2) and (3) of that corollary, with
the notation of \cite[Definition 3.2]{HSY}, express \((a_i)\) and
\((b_i)\) in terms of fixed operators and common isometry systems.
Hence \cite[Theorem 3.8]{HSY} gives \eqref{eq:hsy-hilbert}, with
\[
 \beta=
 \left\|
 \begin{pmatrix}0&d'\\ d&0\end{pmatrix}
 \right\|_{\mathcal C_E(H'\oplus H')},
\]
where \(H'\) is the common Hilbert space in the isometry
representations of \cite[Definition 3.2]{HSY}, and \(d,d'\) are
the operators representing \((a_i)\) and \((b_i)\), respectively.
\end{proof}

For \(E=\ell_{p,q}\), \(q<\infty\), the quasi-norm is order
continuous and has the Fatou property.  Thus \(c_0\not\hookrightarrow
E\), and \cref{prop:hsy} applies; see also
\cite[Proposition A.3]{HSY}.

\section{Disjoint families and finite-dimensional compressions}\label{sec:estimates}

In the Banach range \(1<p<\infty\), the one-sided lower estimate
below follows from Arazy and Lin \cite[Lemma~3.1(ii) and
Corollary~4.3(i)]{ArazyLin}, or from Dodds, Dodds, and Sukochev
\cite[Propositions~5.1(2) and~6.7(1)]{DDS}.  For \(p\leq q\),
the singular-value inequality of \cite[Proposition~4.1]{DDS} also
covers \(0<p\leq1\) and gives constant one for the quasi-norm used here.

\begin{lemma}\label{lem:one-sided}
Let \(0<p<\infty\), \(2\leq q<\infty\), and \(s=\max\{p,q\}\).
There is a constant \(c_{p,q}>0\) such that every finite family
\((x_i)\subset\mathcal C_{p,q}\) with pairwise disjoint left supports
or pairwise disjoint right supports satisfies
\begin{equation}\label{eq:one-sided}
 \left\|\sum_i x_i\right\|_{p,q}
 \geq c_{p,q}\left(\sum_i\|x_i\|_{p,q}^{s}\right)^{1/s}.
\end{equation}
If \(p\leq q\), one may take \(c_{p,q}=1\).
\end{lemma}

\begin{proof}
Suppose the left supports are pairwise disjoint and set \(x=\sum_i x_i\).
Then \(x_i^*x_j=0\) for \(i\ne j\), so \(x^*x=\sum_i x_i^*x_i\).
If \(p\leq q\), apply \cite[Proposition~4.1]{DDS} successively to
these positive summands with the convex function \(\varphi(t)=t^{q/2}\).
For each \(k\geq1\) this gives
\[
 \sum_{j=k}^{\infty}\mu_j(x)^q
 \geq \sum_i\sum_{j=k}^{\infty}\mu_j(x_i)^q.
\]
Put \(w_0=0\) and \(w_k=k^{q/p-1}\) for \(k\geq1\).
Since \(w_k-w_{k-1}\geq0\), multiplying by \(w_k-w_{k-1}\) and
summing over \(k\) yields \(\|x\|_{p,q}^q\geq\sum_i\|x_i\|_{p,q}^q\).

If \(2\leq q<p\), the fully symmetric Banach space
\(L_{p,q}(0,\infty)\) has a lower \(p\)-estimate
\cite[Lemma~4.1(i)]{Randrianantoanina}.  Since
\(L_{p,q}(0,\infty)=(L^2,L^\infty)_{1-2/p,q}\) with equivalent norms,
\cite[Proposition~6.7(1)]{DDS} gives the lower \(p\)-estimate in
\(\mathcal C_{p,q}\).  Taking adjoints gives the assertion for
right-disjoint families in both cases.
\end{proof}

Write \(\kappa\) for the quasi-triangle constant of \(\|\cdot\|_{p,q}\).
For \(0<p,q<\infty\), the closed span of every normalized basic sequence
in \(\mathcal C_{p,q}\) with pairwise disjoint left supports and pairwise
disjoint right supports contains a copy of \(\ell_q\)
\cite[Lemma~4.6]{Randrianantoanina}.  For the commutative result, see
Carothers and Dilworth \cite[Corollary~2.4]{CD}.  For \(q>2\), such a
span cannot be isomorphic to a Hilbert space.  The next lemma then
follows from \cref{prop:hsy} by the standard approximation argument
used by Arazy and Lindenstrauss for the Schatten classes
\cite[Lemma~4, pp.~108--109]{AL}.

\begin{lemma}\label{lem:capture}
Let \(2<q<\infty\) and \(0<p<\infty\).  For each \(K\geq1\) there
is a constant \(\gamma>0\) with the following property.  If
\((x_n)\subset\mathcal C_{p,q}\) is \(\sigma\)-weakly null and
\begin{equation}\label{eq:K-l2}
 K^{-1}\|(t_n)\|_2
 \leq\left\|\sum_nt_nx_n\right\|_{p,q}
 \leq K\|(t_n)\|_2
\end{equation}
for all finitely supported scalar sequences \((t_n)\), then some
\(d\geq1\) satisfies
\begin{equation}\label{eq:capture}
 \liminf_n\max\{\|\pi_dx_n\|_{p,q},\|x_n\pi_d\|_{p,q}\}
 \geq\gamma.
\end{equation}
\end{lemma}

\section{Proof of the classification}\label{sec:main-proof}

Fix \(0<p<\infty\) and \(2<q<\infty\).
Let \((h_j)\) denote the nonconstant \(L_q\)-normalized Haar functions
on \((0,1)\), enumerated level by level.
By a faithful Haar system we mean a system obtained by successively
dividing each atom into two sets of equal measure and taking the
corresponding normalized differences of their indicator functions.
Such a system has the same finite-dimensional distributions as the
standard Haar system.  Hence for any faithful system \((g_j)\) the
map \(h_j\mapsto g_j\) extends to an isometry from the mean-zero
subspace \(L_q^0(0,1)\) onto the closed span of \((g_j)\).

\begin{lemma}\label{lem:annular}
Suppose that \(T:L_q(0,1)\to\mathcal C_{p,q}\) is an isomorphic
embedding.  Then there is another isomorphic embedding
\[
 S:L_q^0(0,1)\longrightarrow\mathcal C_{p,q},
\]
where \(L_q^0\) is the mean-zero hyperplane, and there are integers
\[
 0=M_0<M_1<M_2<\cdots
\]
such that, for the standard normalized Haar basis \((h_j)\) of
\(L_q^0\),
\begin{equation}\label{eq:annular-image}
 Sh_j=z_j
 =
 (\mathscr R_{M_j}-\mathscr R_{M_{j-1}})\,Tg_j
 =
 (\mathscr R_{M_j}-\mathscr R_{M_{j-1}})\,Ah_j.
\end{equation}
Here \((g_j)\) is a faithful Haar system,
\(U:L_q^0(0,1)\to L_q(0,1)\) is the isometric embedding
satisfying \(Uh_j=g_j\), and \(A=TU\).

If
\[
 e_j=\pi_{M_j}-\pi_{M_{j-1}}
\]
and \(d\leq M_{j-1}\), then
\begin{equation}\label{eq:annular-heads}
 \pi_dz_j=\pi_dTg_j\,e_j,\qquad
 z_j\pi_d=e_jTg_j\,\pi_d.
\end{equation}
\end{lemma}

\begin{proof}
Equip \(\mathcal C_{p,q}\) with an equivalent \(\rho\)-subadditive
quasi-norm \(\|\cdot\|_*\leq C\|\cdot\|_{p,q}\), and let
\(a_T>0\) satisfy \(\|Tf\|_*\geq a_T\|f\|_q\).  Writing \(h_j^*\)
for the Haar coordinate functionals, choose \(\eps_j>0\) so that
\(2C^\rho\sum_j\|h_j^*\|^\rho\eps_j^\rho<a_T^\rho/2\).
At step \(j\), let \(D\) be the parent atom on which \(g_j\) is to
be supported; for \(j=1\), take \(D=(0,1)\).  On \(D\) with
normalized measure choose an independent Rademacher sequence.
Multiplying each sign by \(|D|^{-1/q}\) and extending by zero off
\(D\) gives a weakly null sequence of \(L_q\)-unit vectors.  The range of
\(\mathscr R_{M_{j-1}}T\) is finite-dimensional.  We may therefore
choose a normalized sign \(g_j\) so that
\begin{equation}\label{eq:old-square-small}
 \|\mathscr R_{M_{j-1}}Tg_j\|_{p,q}<\eps_j.
\end{equation}
The positive and negative sets of \(g_j\) partition \(D\) into two
child atoms of equal measure.
By \eqref{eq:square-compression}, choose \(M_j>M_{j-1}\) so that
\begin{equation}\label{eq:new-square-tail}
 \|(I-\mathscr R_{M_j})Tg_j\|_{p,q}<\eps_j.
\end{equation}
Continuing inductively gives a faithful Haar system \((g_j)\).
Let \(U\) be the isometric embedding determined by \(Uh_j=g_j\),
and put \(A=TU\).  The difference \(Ah_j-z_j\) is the sum of the
two errors in \eqref{eq:old-square-small} and \eqref{eq:new-square-tail}, so
\[
 \|Ah_j-z_j\|_*^\rho\leq 2C^\rho\eps_j^\rho .
\]
Thus the perturbation \(Df=\sum_jh_j^*(f)(Ah_j-z_j)\) satisfies
\[
 \|Df\|_*^\rho
 \leq 2C^\rho\|f\|_q^\rho\sum_j\|h_j^*\|^\rho\eps_j^\rho
 \leq (a_T^\rho/2)\|f\|_q^\rho .
\]
Since \(U\) is isometric, \(\|Af\|_*\geq a_T\|f\|_q\).
It follows that \(S=A-D\) satisfies
\(\|Sf\|_*^\rho\geq\|Af\|_*^\rho-\|Df\|_*^\rho
\geq(a_T^\rho/2)\|f\|_q^\rho\), and \(Sh_j=z_j\).

Finally, if \(d\leq M_{j-1}\), direct multiplication in
\eqref{eq:annular-image} gives
\[
 \pi_dz_j
 =\pi_dTg_j(\pi_{M_j}-\pi_{M_{j-1}})
 =\pi_dTg_j e_j,
\]
and similarly \(z_j\pi_d=e_jTg_j\pi_d\).
\end{proof}

\begin{remark}\label{rem:late-level}
Fix \(d\geq0\).  At every sufficiently large Haar level, all indices
\(\ell\) satisfy \(M_{\ell-1}\geq d\).  If
\(h=\sum_{\ell\in F}c_\ell h_\ell\) belongs to such a level, put
\(e_F=\sum_{\ell\in F}e_\ell\).  By \eqref{eq:annular-heads},
\(\pi_dSh=(\pi_dSh)e_F\) and \(Sh\pi_d=e_F(Sh\pi_d)\).
Thus disjoint index sets \(F\) give orthogonal right supports
after left compression and orthogonal left supports after right compression.
\end{remark}

We next prove a selection lemma for symmetric operator ideals.
The constants in this lemma are independent of the underlying
sequence space.

\begin{lemma}\label{lem:selection}
Let \(N\geq1\) and \(d\geq0\) be integers.  Let
\(x_1,\ldots,x_N\in\mathcal C_F\), where \(F\) is any symmetric
quasi-Banach sequence space.  Suppose that the operators
\[
 u_j=\pi_dx_j
\]
have pairwise orthogonal right supports, and the operators
\[
 v_j=x_j\pi_d
\]
have pairwise orthogonal left supports.  There are numerical constants
\(C_0>0\) and \(C_1>0\) with the following property.  For each
positive integer \(m\) satisfying
\[
 C_0d\log(ed)\leq m\leq N/4,
\]
with the lower bound omitted when \(d=0\), there exists
\(A\subset\{1,\ldots,N\}\), \(|A|=m\), such that
\begin{align}
 \left\|\sum_{j\in A}u_j\right\|_{\mathcal C_F}
 &\leq C_1\sqrt{\frac mN}
 \left\|\sum_{j=1}^Nu_j\right\|_{\mathcal C_F},
 \label{eq:select-row}\\
 \left\|\sum_{j\in A}v_j\right\|_{\mathcal C_F}
 &\leq C_1\sqrt{\frac mN}
 \left\|\sum_{j=1}^Nv_j\right\|_{\mathcal C_F}.
 \label{eq:select-column}
\end{align}
\end{lemma}

\begin{proof}
If \(d=0\), both conclusions are immediate; hence assume \(d\geq1\).
Put
\[
 A_j^{(r)}=u_ju_j^*,\quad S_r=\sum_jA_j^{(r)},\qquad
 A_j^{(c)}=v_j^*v_j,\quad S_c=\sum_jA_j^{(c)}.
\]
On \(\operatorname{ran}S_r\oplus\operatorname{ran}S_c\), a space of
dimension \(r=\operatorname{rank}S_r+\operatorname{rank}S_c\leq2d\), define
\[
 B_j=
 S_r^{-1/2}A_j^{(r)}S_r^{-1/2}
 \ \oplus\
 S_c^{-1/2}A_j^{(c)}S_c^{-1/2}.
\]
The generalized inverses are taken on the respective supports.
Then \(B_j\geq0\), \(\sum_jB_j=I_r\), and \(r\leq2d\).
If \(r=0\), every \(u_j\) and \(v_j\) is zero and the conclusion
holds for any \(m\)-element set.  Assume \(r\geq1\), and put
\[
 G=\{j:\|B_j\|_\infty\leq2r/N\}.
\]
Then \(G\) has at least \(N/2\) elements, since
\(\sum_j\Tr(B_j)=r\).

Select the indices in \(G\) independently with probability
\(\rho=2m/|G|\leq1\).  The matrix Bernstein inequality and the scalar
Chernoff bound imply that, for a sufficiently large numerical
constant \(C\), the condition \(m\geq Cr\log(er)\) gives a
realization with at least \(m\) selected indices and
\[
 \sum_{\text{selected }j}B_j\leq C(m/N)I_r.
\]
Indeed, the summands are bounded by \(2r/N\), the variance is
bounded by a numerical multiple of \(mr/N^2\), and the probability
that either conclusion fails is at most
\(2r\exp(-cm/r)+\exp(-c'm)\); see
\cite[Theorem 6.1(ii)]{Tropp}.
Since \(r\leq2d\), the hypothesis \(m\geq C_0d\log(ed)\) suffices
when \(C_0\) is sufficiently large.
Choose any \(m\) of the selected indices and denote their set by
\(A\).  Since \(B_j\geq0\), we have \(\sum_{j\in A}B_j\leq C(m/N)I_r\).

Conjugating the first diagonal block of this inequality by \(S_r^{1/2}\)
and the second by \(S_c^{1/2}\) gives
\[
 \sum_{j\in A}u_ju_j^*
 \leq C(m/N)\sum_ju_ju_j^*,
\quad
 \sum_{j\in A}v_j^*v_j
 \leq C(m/N)\sum_jv_j^*v_j.
\]
Since the \(u_j\)'s have orthogonal right supports,
\[
 \left(\sum_{j\in A}u_j\right)
 \left(\sum_{j\in A}u_j\right)^*
 =\sum_{j\in A}u_ju_j^*,
\]
and the corresponding identity for \(v_j^*v_j\) follows from
orthogonality of the left supports of the \(v_j\)'s.
The inequalities \eqref{eq:select-row}--\eqref{eq:select-column}
now follow from eigenvalue monotonicity and symmetry of the
quasi-norm.
\end{proof}

\begin{proof}[Proof of the case \(q>2\)]
Assume that \(T:L_q(0,1)\to\mathcal C_{p,q}\) is an isomorphic
embedding.  By \cref{lem:annular}, there is an isomorphic embedding
\(S:L_q^0(0,1)\to\mathcal C_{p,q}\) satisfying
\eqref{eq:annular-image}--\eqref{eq:annular-heads}.  Fix
\(0<a\leq b<\infty\) with
\begin{equation}\label{eq:S-bounds}
 a\|f\|_q\leq\|Sf\|_{p,q}\leq b\|f\|_q.
\end{equation}

Let \(\mathcal F_k\) be the sigma-algebra generated by the dyadic
intervals of length \(2^{-k}\), and let \(r_n\) be the sign that
equals \(+1\) and \(-1\) on the left and right halves, respectively,
of each dyadic interval of length \(2^{-(n-1)}\).  For
\(B\in\mathcal F_k\) with \(|B|>0\) and \(n>k\), put
\[
 s_n^B=|B|^{-1/q}r_n1_B.
\]
The restrictions \(r_n|_B\) are independent Rademacher variables
on \((B,|B|^{-1}dt)\), and \(\|s_n^B\|_{L_q(0,1)}=1\).
By Khintchine's inequality and
\eqref{eq:S-bounds}, the sequence \((S(s_n^B))_{n>k}\) is
\(K\)-equivalent to the unit vector basis of \(\ell_2\), with
\(K\) depending only on \(a,b,q\), not on \(k\) or \(B\).
Fix \(\gamma>0\) given by \cref{lem:capture} for this \(K\).

Choose a dyadic number \(t\in(0,1/4)\), fixed throughout the proof,
such that
\begin{equation}\label{eq:t-small}
 C_1b\,t^{1/2-1/q}\leq\frac{\gamma}{4\kappa}.
\end{equation}
Here \(C_1\) is the constant in \cref{lem:selection}.  Set
\(\delta=t^{1/q}\gamma/(4\kappa)\).

We construct inductively integers
\[
 0=d_0<d_1<\cdots<d_s<\cdots,
\]
nested infinite index sets \(\mathbb N=\Lambda_0\supseteq\Lambda_1
\supseteq\cdots\), and a choice of left or right compression at
each step.  For each \(1\leq i\leq s\), one of the following
inequalities must hold for every \(n\in\Lambda_s\):
\begin{equation}\label{eq:recorded-shell}
 \|( \pi_{d_i}-\pi_{d_{i-1}})S(r_n)\|_{p,q}\geq\delta,
\end{equation}
or
\begin{equation}\label{eq:recorded-column-shell}
 \|S(r_n)( \pi_{d_i}-\pi_{d_{i-1}})\|_{p,q}\geq\delta.
\end{equation}
The side of compression at step \(i\)
is independent of \(n\) and is retained at all subsequent steps.

Given \(d=d_s\) and \(\Lambda_s\), choose \(k\) sufficiently large that
\(t2^k\) is an integer
and, with \(N=2^k\),
\[
 tN\geq C_0d\log(ed).
\]
When \(d=0\), the lower bound is omitted.
Let \(I_{k,j}\), \(1\leq j\leq N\), be the dyadic intervals of
length \(N^{-1}\).  For \(n\in\Lambda_s\) with \(n>k\), put
\[
 \widehat r_{n,k,j}=N^{1/q}r_n1_{I_{k,j}}.
\]
Thus
\[
 r_n=N^{-1/q}\sum_{j=1}^N\widehat r_{n,k,j}.
\]
Each \(r_n1_{I_{k,j}}\) is a linear combination of Haar functions
supported on the intervals of length \(2^{-(n-1)}\) contained in
\(I_{k,j}\).  Distinct \(j\)'s therefore use disjoint sets of Haar
indices.  Discard finitely many \(n\) so that \(M_{\ell-1}\geq d\)
for every Haar index \(\ell\) at level \(n-1\).  By \cref{rem:late-level}, for
\[
 x_{n,j}=S(\widehat r_{n,k,j})
\]
the left compressions \(\pi_dx_{n,j}\) are right-disjoint and the
right compressions \(x_{n,j}\pi_d\) are left-disjoint.  Apply
\cref{lem:selection} with \(m=tN\).  For every remaining \(n\) there is
\(A_n\subset\{1,\ldots,N\}\), \(|A_n|=tN\), satisfying both
selection estimates.  Since there are only finitely many possible
sets, pass to an infinite subsequence on which \(A_n=A\) is fixed.
Put
\[
 B=\bigcup_{j\in A}I_{k,j},\qquad
 s_n=t^{-1/q}r_n1_B.
\]
Multiplying the estimates in \cref{lem:selection} by
\(t^{-1/q}N^{-1/q}\) and using \eqref{eq:S-bounds}, we obtain
\begin{equation}\label{eq:two-old-heads}
 \|\pi_dS(s_n)\|_{p,q},\quad
 \|S(s_n)\pi_d\|_{p,q}
 \leq C_1b\,t^{1/2-1/q}.
\end{equation}

The sequence \((s_n)\) is weakly null in \(L_q\).
Every matrix coefficient is a continuous functional on
\(\mathcal C_{p,q}\), because
\(\|x\|_\infty\lesssim\|x\|_{p,q}\).  Hence
\((S(s_n))\) converges to zero in the weak-operator topology.
Moreover, \(\sup_n\|S(s_n)\|_\infty\leq Cb\); on an
operator-norm bounded set, weak-operator convergence implies
\(\sigma\)-weak convergence.  Thus \((S(s_n))\) is
\(\sigma\)-weakly null.  By
\cref{lem:capture}, there is an integer \(d'\) for which the
maximum of the two compression norms has liminf at least \(\gamma\).
After passage to an infinite subsequence, one side is fixed: either
\[
 \|\pi_{d'}S(s_n)\|_{p,q}\geq\gamma/2
\]
for all remaining \(n\), or
\[
 \|S(s_n)\pi_{d'}\|_{p,q}\geq\gamma/2
\]
for all remaining \(n\).  By \eqref{eq:t-small} and
\eqref{eq:two-old-heads}, both compressions at \(d\) have norm at
most \(\gamma/(4\kappa)\).  Since the projections \(\pi_d\)
increase with \(d\), we must have \(d'>d\).
The quasi-triangle inequality gives, in the first case,
\[
 \|(\pi_{d'}-\pi_d)S(s_n)\|_{p,q}
 \geq\frac{\gamma}{4\kappa},
\]
and in the second case,
\[
 \|S(s_n)(\pi_{d'}-\pi_d)\|_{p,q}
 \geq\frac{\gamma}{4\kappa}.
\]

Discard finitely many more \(n\) so that \(M_{\ell-1}\geq d'\) for
every Haar index \(\ell\) at level \(n-1\).  Let \(\Lambda_{s+1}\)
consist of the remaining indices \(n\), and put \(Q=\pi_{d'}-\pi_d\).
Since
\[
 r_n=t^{1/q}s_n+r_n1_{B^c},
\]
\eqref{eq:annular-heads} gives orthogonal right supports for
\(QS(s_n)\) and \(QS(r_n1_{B^c})\), and orthogonal left supports for
\(S(s_n)Q\) and \(S(r_n1_{B^c})Q\).  For a left compression,
singular-value monotonicity gives
\[
 \|QS(r_n)\|_{p,q}\geq t^{1/q}\|QS(s_n)\|_{p,q}\geq\delta .
\]
For a right compression, the corresponding bound is
\[
 \|S(r_n)Q\|_{p,q}\geq t^{1/q}\|S(s_n)Q\|_{p,q}\geq\delta .
\]
Taking \(d_{s+1}=d'\) completes the induction.

After \(L\) steps, choose \(n\in\Lambda_L\).  At least \(L/2\)
of the first \(L\) steps use the same side of compression.
Suppose first that these are left compressions, and let \(J\)
be the set of their indices.
The corresponding operators
\[
 (\pi_{d_i}-\pi_{d_{i-1}})S(r_n)
\]
are left-disjoint.  Their sum is a left compression of \(S(r_n)\).
Set \(s=\max\{p,q\}\).  Applying \cref{lem:one-sided} to the family
indexed by \(J\) gives
\[
 b\geq\|S(r_n)\|_{p,q}
 \geq
 \left\|\sum_{i\in J}
 (\pi_{d_i}-\pi_{d_{i-1}})S(r_n)\right\|_{p,q}
 \geq c_{p,q}\delta(L/2)^{1/s}.
\]
The case of right compressions follows from the right-disjoint
estimate in the same way.  The index \(n\) may depend on \(L\),
but \(b,\delta,s\) and \(c_{p,q}\) do not.  Letting
\(L\to\infty\) contradicts the upper bound \(b\).
\end{proof}

The converse in \cref{thm:classification} follows from the isometric
identification of \(L_2(0,1)\) with \(\ell_2\): a matrix supported on
the first row has ideal quasi-norm equal to the \(\ell_2\)-norm of that row.


\begin{thebibliography}{99}

\bibitem{ArazyLin}
J.~Arazy and P.-K.~Lin,
\emph{On \(p\)-convexity and \(q\)-concavity of unitary matrix spaces},
Integral Equations Operator Theory \textbf{8} (1985), no.~3, 295--313,
\href{https://doi.org/10.1007/BF01202902}{doi:10.1007/BF01202902}.

\bibitem{AL}
J.~Arazy and J.~Lindenstrauss,
\emph{Some linear topological properties of the spaces
\(\mathcal C_p\) of operators on Hilbert space},
Compositio Math. \textbf{30} (1975), 81--111.

\bibitem{AHS}
S.~Astashkin, J.~Huang, and F.~Sukochev,
\emph{Lack of isomorphic embeddings of symmetric function spaces into
operator ideals},
J. Funct. Anal. \textbf{280} (2021), no.~5, Paper No.~108895.

\bibitem{CD}
N.~L. Carothers and S.~J. Dilworth,
\emph{Subspaces of \(L_{p,q}\)},
Proc. Amer. Math. Soc. \textbf{104} (1988), no.~2, 537--545,
\href{https://doi.org/10.1090/S0002-9939-1988-0962825-8}
{doi:10.1090/S0002-9939-1988-0962825-8}.

\bibitem{DDS}
P.~G. Dodds, T.~K. Dodds, and F.~A. Sukochev,
\emph{On \(p\)-convexity and \(q\)-concavity in non-commutative
symmetric spaces},
Integral Equations Operator Theory \textbf{78} (2014), no.~1, 91--114,
\href{https://doi.org/10.1007/s00020-013-2082-0}
{doi:10.1007/s00020-013-2082-0}.

\bibitem{HSY}
J.~Huang, F.~Sukochev, and Z.~Yu,
\emph{Arazy-type decomposition theorem for bounded linear operators
and commutators on the trace class},
arXiv:2602.09579v1 [math.FA], 10 February 2026, preprint.

\bibitem{HSZ}
J.~Huang, F.~Sukochev, and D.~Zanin,
\emph{Lack of isomorphic embedding from Lorentz function spaces into
Lorentz operator ideals},
Adv. Math. \textbf{504} (2026), Paper No.~111276,
\href{https://doi.org/10.1016/j.aim.2026.111276}
{doi:10.1016/j.aim.2026.111276}.

\bibitem{Randrianantoanina}
N.~Randrianantoanina,
\emph{Embeddings of \(\ell_p\) into non-commutative spaces},
J. Aust. Math. Soc. \textbf{74} (2003), no.~3, 331--350,
\href{https://doi.org/10.1017/S1446788700003359}
{doi:10.1017/S1446788700003359}.

\bibitem{Tropp}
J.~A.~Tropp,
\emph{User-friendly tail bounds for sums of random matrices},
Found. Comput. Math. \textbf{12} (2012), no.~4, 389--434.

\end{thebibliography}
\end{document}